\documentclass{article}

\usepackage[french,british]{babel}

\usepackage[T1]{fontenc}

\usepackage{amsmath}
\DeclareMathOperator{\tr}{tr}
\DeclareMathOperator{\RE}{Re}

\usepackage{amsfonts}
\usepackage{amsbsy}
\usepackage{amssymb}

\usepackage{amsthm}

\theoremstyle{plain}

\theoremstyle{definition}
\newtheorem{definition}{Definition}

\newtheorem{exm}{Example}
\newtheorem{lemma}{Lemma}

\newtheorem{Corollary}{Corollary}

\newtheorem{thm}{Theorem}

\begin{document}

\title{Minimax Extrapolation Problem For Harmonizable Stable Sequences With Noise Observations}

\author{Mikhail Moklyachuk$^1$
\footnote{Corresponding email:Moklyachuk@gmail.com%
 \newline 1 Department of Probability Theory, Statistics and Actuarial
Mathematics, Taras Shevchenko National University of Kyiv, Kyiv, Ukraine%
%
 }
             and Vitalii Ostapenko$^{1}$
 }

\date{\today}

\maketitle

\renewcommand{\abstractname}{Abstract}
\begin{abstract}
 We consider the problem of optimal linear estimation of the functional $A \xi~=~\sum_{j = 0}^{\infty} a_j \xi_j$ that depends on the
unknown values $\xi_j,j=0,1,\dots, $ of a  random sequence $\{\xi_j,j\in\mathbb Z\}$ from observations of the sequence
 $\{\xi_k+\eta_k,k\in\mathbb Z\}$ at points
$k = -1, -2, \dots$, where $\{\xi_k,k\in\mathbb Z\}$ and $\{\eta_k,k\in\mathbb Z\}$ are mutually independent harmonizable symmetric $\alpha$-stable random sequences
which have the spectral densities $f(\theta)>0$ and $g(\theta)>0$ satisfying the minimality condition.
 The problem is investigated under the condition of spectral certainty as well as under the condition of spectral uncertainty.
 Formulas for calculating the value of the error and the spectral characteristic of the optimal linear
estimate of the functional are derived under the condition of spectral certainty where spectral densities of the sequences are exactly known.
 In the case of spectral uncertainty where spectral densities of
the sequences are not exactly known, while a class of admissible spectral densities
is given, relations that determine the least favorable spectral densities and the minimax spectral
characteristic are derived.
\end{abstract}

\vspace{2ex}

\vspace{2ex}

\selectlanguage{british}

\section{Introduction}

 The classical methods of finding solutions to extrapolation, interpolation and filtering problems for stationary stochastic processes and sequences were developed by Kolmogorov (see selected works of Kolmogorov (1992)), Wiener (see the book by Wiener (1966)), Yaglom (see, for example, books by Yaglom (1987a, 1987b).
The problem of estimation of the unknown values of harmonizable random sequences and processes were investigated in papers by Cambanis (1983),
 Cambanis and Soltani (1984), Hosoya (1982).
The interpolation problem for harmonizable symmetric $\alpha$-stable random sequences were investigated in papers by
Weron (1985) and Pourahmadi (1984).
Most of results concerning estimation of the unknown (missed) values of stochastic processes are based on the assumption that spectral densities of processes are exactly known. In practice, however, complete information on the spectral densities is impossible in most cases.
In such situations one finds parametric or nonparametric estimates of the unknown spectral densities. Then the classical estimation method is applied under the assumption that the estimated densities are true. This procedure can result in significant increasing of the value of error as Vastola and Poor (1983)  have demonstrated with the help of some examples.
 This is a reason to search estimates which are optimal for all densities from a certain class of admissible spectral densities.
 These estimates are called minimax since they minimize the maximal value of the error.
 A survey of results in minimax (robust) methods of data processing can be found in the paper by Kassam and Poor (1985).
 The paper by Grenander (1957) should be marked as the first one where the minimax extrapolation problem for stationary processes was formulated and solved.
Later Franke and Poor (Franke and Poor, 1984; Franke, 1985) applied the convex optimization methods for investigation the minimax-robust extrapolation and interpolation problems.
In papers by Moklyachuk (1990 -- 2008a) of the minimax-robust extrapolation, interpolation and filtering  problems are studied for stationary processes. The papers by Moklyachuk and Masyutka
(2006 -- 2012) are dedicated to minimax-robust extrapolation, interpolation and filtering  problems for vector-valued stationary processes and sequences.
 Dubovetska et al. (2012) solved the problem of minimax-robust interpolation for another generalization of stationary processes --  periodically correlated sequences. In the papers by  Dubovetska and Moklyachuk (2013 -- 2014),  Moklyachuk and Golichenko (2016) the minimax-robust extrapolation, interpolation and filtering  problems for periodically correlated processes are investigated.
 The minimax-robust extrapolation, interpolation and filtering  problems for stochastic  sequences and random processes with $n$th stationary increments are investigated by Luz and Moklyachuk (Luz and Moklyachuk, 2012 -- 2016; Moklyachuk and Luz, 2013).

In this paper the problem of optimal estimation is investigated for the linear functional $A \xi~=~\sum_{j = 0}^{\infty} a_j \xi_j$ that depends on the
unknown values of a random sequence $\{\xi_j,j\in\mathbb Z\}$ based on observations of the sequence
 $\{\xi_k+\eta_k,k\in\mathbb Z\}$ at points
$k = -1,  -2, \dots$, where $\{\xi_k,k\in\mathbb Z\}$ and $\{\eta_k,k\in\mathbb Z\}$ are mutually independent harmonizable symmetric $\alpha$-stable random sequences
which have the spectral densities $f(\theta)>0$ and $g(\theta)>0$ satisfying the minimality condition.
The problem is investigated under the condition of spectral certainty as well as under the condition of spectral uncertainty.
 Formulas for calculating the value of the error and the spectral characteristic of the optimal linear
estimate of the functional are derived under the condition of spectral certainty where spectral densities of the sequences are exactly known.
 In the case of spectral uncertainty where spectral densities of
the sequences are not exactly known while sets of admissible spectral densities are available, relations which determine the least favorable densities and the minimax-robust spectral characteristics for different classes of spectral densities are derived.

\section{Harmonizable symmetric $\alpha$-stable random sequences. Basic properties}\label{sec:fields}

 \begin{definition}[symmetric $\alpha$-stable random variable]
 A real random variable $\xi$ is said to be symmetric $\alpha$-stable, $S\alpha S$, if its characteristic function has the form $E exp(it\xi) = exp(-c|t|^{\alpha})$ for some $c \geq 0$ and $0 < \alpha \leq 2.$
 The real random variables $\xi_1,\xi_2,\dots,\xi_n$ are jointly $S\alpha S$ if all linear combinations
 $\sum_{k=1}^{n}a_k\xi_k$ are $S\alpha S$, or ,equivalently, if the characteristic function of $\vec{\xi}=(\xi_1,\dots,\xi_n)$ is of the form
 $\phi_{\vec{\xi}}(\vec{t}) = E exp(i \sum t_k \xi_k) = exp\{-\int|\sum t_k x_k|^{\alpha}d \Gamma_{\vec{\xi}}(\vec{x})\},$ where $t_1,\dots, t_n$ are real numbers and $\Gamma_{\xi}$ is a symmetric measure defined on the unit sphere $S_n \in R^n$ (Cambanis, 1983).
 \end{definition}

 \begin{definition}[symmetric $\alpha$-stable stochastic sequence]
A stochastic sequence $\{\xi_k,k\in\mathbb Z\}$ is called symmetric $\alpha$-stable, $S\alpha S$, stochastic sequence, if all  linear combinations $\sum_{k=1}^{n}a_k\xi_k$ are $S\alpha S$ random variables.
 \end{definition}

 For jointly $S{\alpha}S$ random variables $\xi= \xi_1 + i\xi_2$ and $\eta= \eta_1 + i \eta_2$ the covariation of $\xi$ with $\eta$ is defined as (Cambanis, 1983)
 $$[\xi,\eta]_{\alpha} = \int_{S_4} (x_1 + i x_2)(y_1 + i y_2)^{<\alpha-1>} d \Gamma_{\xi_1,\xi_2,\eta_1,\eta_2}(x_1, x_2, y_1, y_2),$$

 \noindent where $z^{<\beta>} = |z|^{\beta - 1} \bar{z} $ for a complex number $z$ and $\beta > 0.$
 The covariation in general is not symmetric and linear on second argument (Weron, 1985). For $\xi, \xi_1, \xi_2, \eta$ jointly $S\alpha S$ we have
 $$[\xi_1 + \xi_2, \eta]_{\alpha} = [\xi_1, \eta]_{\alpha} + [\xi_2, \eta]_{\alpha},$$

 \begin{equation}\label{eq:inequality}
 |[\xi, \eta]_{\alpha}| \leq ||\xi||_{\alpha} ||\eta||_{\alpha}^{\alpha - 1}
 \end{equation}
 and $||\xi||_{\alpha} = [\xi, \xi]_{\alpha}^{1/\alpha}$ is a norm in a linear space of $S\alpha S$ random variables which is equivalent to convergence in probability.
 It should be noted that $||\cdot||_\alpha$ is not necessarily the usual $L^{\alpha}$ norm.

 Here is the simplest properties of the function $z^{<\beta>}.$
\begin{lemma}
 Let $z, x, y $ be complex numbers, $\beta > 0$. Then the following properties hold true:
 \begin{itemize}
 \item $|z|^{<\beta>} = z \cdot z^{<\beta - 1>},$
 \item $\left||z|^{<\beta>}\right| = \left|z\right|^{<\beta>},$
 \item if $z^{<\beta>} = v$, thet $z = v^{<1/\beta>} = |v|^{(1-\beta)/\beta}\bar{v},$
 \item $z^{<1>} = \bar{z},$
 \item if $z \neq 0$, then $z^{<\alpha>} z^{<\beta>} = \frac{\bar{z}}{|z|} z^{<\alpha + \beta>},$
 \item if $z \neq 0$, then $\frac{z^{<\alpha>}}{z^{<\beta>}} = \frac{z}{|z|} z^{<\alpha - \beta>},$
 \item $(c z)^{<\alpha>} = c^{\alpha} z^{<\alpha>}, c \in \mathbb{R},$
 \item $(z^{<\alpha>})^{<\beta>} = {\bar{z}}^{<\alpha \beta>},$
 \item $(xy)^{<\alpha>} = x^{<\alpha>} y^{<\alpha>},$
 \item $(z^{\alpha})^{<\beta>} = (z^{<\beta>})^{\alpha},$
 \item $(z^{<\alpha>})^{\beta} = (z^{\beta})^{<\alpha>},$
 \item $|z^{<\alpha>}|^{\beta} = |z|^{\alpha \beta},$
 \item $(x + y)^{<\alpha>} = \bar{x}|x + y|^{\alpha - 1} + \bar{y}|x + y|^{\alpha - 1}.$
\end{itemize}
\end{lemma}

Let $Z =\{Z(t): -\infty < t < \infty\}$ be a complex $S{\alpha}S$ process with independent increments. The spectral measure of the process $Z$ is defined as $\mu\{(s, t]\} =\|Z(t) - Z(s)\|_{\alpha}^{\alpha}.$ The integrals $\int f(t)dZ(t)$ can be defined for all $f \in L^{\alpha}(\mu)$ with properties (see Cambanis, 1983; Cambanis and Soltani, 1984; Hosoya, 1982):
\[
 \left\|\int f(t) d Z(t)\right\|^{\alpha}_{\alpha} = \int |f(t)|^{\alpha}d \mu,
\]
\begin{equation}\label{eq:norm_equality}
 \left[\int f(t) d Z(t), \int g(t) d Z(t)\right]_\alpha = \int f (t) (g(t))^{<\alpha - 1>} d \mu.
\end{equation}

\begin{definition}[Harmonizable symmetric $\alpha$-stable stochastic sequence]
 A $S\alpha S$ stochastic sequence $\{\xi_n,n\in\mathbb Z\}$ is said to be harmonizable, $HS{\alpha}S$, if there exists a $S\alpha S$ process $Z = \{Z(\theta); \theta \in [-\pi, \pi]\}$ with independent increments and finite spectral measure $\mu$ such that sequence $\xi_n$ has the spectral representation
 $$\xi_n = \int_{-\pi}^{\pi}e^{in\theta}dZ(\theta), \quad n \in \mathbb{Z},$$

 \noindent and the covariation has the representation
 $$[\xi_n, \xi_m]_\alpha = \int_{-\pi}^{\pi}e^{i(n-m)\theta}d \mu(\theta), \quad m, n \in \mathbb{Z}.$$
\end{definition}

Note that a $HS\alpha S$ stochastic sequence is not necessarily stationary even second order stationary, but for $\alpha = 2$ the $HS\alpha S$ sequences are stationary with Gaussian distribution.
In this article we consider the case where $1<\alpha \leq 2$.

Denote by $H(\xi)$ the time domain of the $HS\alpha S$ sequence $\{\xi_n,n\in\mathbb Z\}$, which is a closed in the norm $\|\cdot\|_{\alpha}$ linear manifold generated by all values of the $HS\alpha S$ sequence $\{\xi_n,n\in\mathbb Z\}$.
It follows from the spectral representation of the $HS\alpha S$ sequence $\{\xi_n,n\in\mathbb Z\}$
that the mapping $\xi_n\leftrightarrow e^{in\theta},n\in\mathbb Z,$ extents to an isomorphism between the spaces $H(\xi)$ and $L^{\alpha}(\mu)$. Under this isomorphism to each $\eta \in H(\xi)$ corresponds a unique $f\in L^{\alpha}(\mu)$ such that $\eta=\int_{-\pi}^{\pi}f(\theta)dZ(\theta)$.

For a closed linear subspace $M \subseteq L^{\alpha}(\mu)$ and $f \in L^{\alpha}(\mu)$, there exists a unique element from $M$ which minimizes the distance to $f$. This element is called projection of $f$ onto $M$ or the best approximation of $f$ in $M$. This projection is denoted by $P_M f$ and is uniquely determined by the condition (Singer, 1970)
\begin{equation}
 \int_{-\pi}^{\pi} g \left(f - P _M f\right)^{<\alpha - 1>}d \mu = 0,\quad g \in M.
\end{equation}

Similarly, for $HS \alpha S$ stochastic sequence $\{\xi_n,n\in\mathbb Z\}$ and a closed linear subspace $H^-(\xi)$ of the space $H(\xi)$ there is a uniquely determined element $\hat{\xi}_n \in H^-(\xi)$ which minimizes the distance to $\xi_n$ and is uniquely determined from the condition
 \begin{equation}\label{eq:ortogonal}
 \left[\eta, \xi_n - \hat{\xi}_n\right]_{\alpha} = 0,\quad \eta \in H^-(\xi).
 \end{equation}

From linearity of the covariation with respect to the first argument from this relation we have that
 \begin{equation}\label{eq:linearity}
 ||\xi_n - \hat{\xi}_n||_{\alpha}^{\alpha} = \left[\xi_n,\xi_n-\hat{\xi}_n\right]_{\alpha}-\left[\hat{\xi}_n, \xi_n-\hat{\xi}_n\right]_{\alpha} = \left[\xi_n, \xi_n-\hat{\xi}_n\right]_{\alpha}.
 \end{equation}

This relation plays a fundamental role in the characterization of minimal $HS \alpha S$ stochastic sequences $\{\xi_n,n\in\mathbb Z\}$.

\section{Extrapolation problem. Projection approach}

Consider the problem of the optimal estimation of the linear functional
$$A \xi~=~\sum_{j = 0}^{\infty} a_j \xi_j = \int_{-\pi}^{\pi} A(e^{i\theta})dZ^{\xi}(\theta),$$
$$A(e^{i\theta})=\sum_{j = 0}^{\infty}a_j e^{ij\theta},
$$
that depends on the
unknown values $\xi_j,j=0,1,\dots,$ of a harmonizable symmetric $\alpha$-stable random sequence $\{\xi_j,j\in\mathbb Z\}$ from observations of the sequence $\{\xi_k+\eta_k,k\in\mathbb Z\}$
 at points
$k = -1,  -2, \dots$.

We will suppose that the sequence $ \{ {a}_j: j=0,1, \ldots \}$ which determines the functional $A  { \xi} $ satisfies conditions
 \begin{equation} \label{cond-on-aj}
  \sum_{j=0}^{ \infty} \left|a_j \right|  < \infty , \quad \sum_{j=0}^{ \infty}(j+1) \left | {a}_j \right | ^{2}  < \infty.
 \end{equation}
 The first condition ensures that the functional  $A\xi$ has a finite second moment. The second condition ensures the compactness in  $\ell_2$ of the operators that will be defined below.

We consider the problem for mutually independent harmonizable symmetric $\alpha$-stable random sequences $\{\xi_k,k\in\mathbb Z\}$ and $\{\eta_k,k\in\mathbb Z\}$ which have
absolutely continuous spectral measures and the spectral densities $f(\theta)>0$ and $g(\theta)>0$ satisfying the minimality condition (Kolmogorov, 1992; Rozanov, 1967; Salehi, 1979; Pourahmadi, 1984; Weron, 1985)
\begin{equation}\label{eq:minimality}
 \int_{-\pi}^{\pi} (f(\theta)+g(\theta))^{-1/(\alpha-1)}d\theta<\infty.
\end{equation}

Denote by $H^{-}(\xi+\eta)$ the closed in the $||\cdot||_\alpha$ norm linear manifold generated by values of the harmonizable symmetric $\alpha$-stable random sequence $\xi_k+\eta_k,
k = -1,  -2, \dots$ in the space $H(\xi+\eta)$ generated by all values of the $HS \alpha S$ sequence $\{\xi_k+\eta_k,k\in\mathbb Z\}$.

The optimal estimate $\hat{A} \xi$ of the functional ${A} \xi$ is a projection of ${A} \xi$ on the manifold $H^{-}(\xi+\eta)$. It is determined by relations
 $$[\zeta, A \xi - \hat{A} \xi]_{\alpha} = 0, \quad \forall \zeta \in H^{-}(\xi+\eta),$$
or, equivalently, by relations
 \begin{equation}\label{ortog1}
 [\xi_k+\eta_k, A \xi - \hat{A} \xi]_{\alpha} = 0, \quad \forall k = -1,  -2, \dots.
 \end{equation}

It follows from the isomorphism between the spaces $H(\xi+\eta)$ and $L^{\alpha}(f+g)$ that the optimal estimate $\hat{A} \xi$ of the functional ${A} \xi$ is of the form
 \begin{equation}\label{estim2}
 \hat{A} \xi = \int_{-\pi}^{\pi} {h}(\theta) \left(d Z^{\xi}(\theta) + dZ^{\eta}(\theta) \right).
 \end{equation}
 It is  determined by the spectral characteristic ${h}(\theta)$ of the estimate which is from the subspace $L^{\alpha}_{-}(f+g)$ of the $L^{\alpha}(f+g)$ space generated by functions
 $e^{ik\theta},k = -1,  -2, \dots$. The spectral characteristic ${h}(\theta)$ of the optimal estimate satisfies the
 following equations
 \begin{equation}\label{ortog2}
 \int_{-\pi}^{\pi} e^{i\theta k} \left[\left( A(e^{i\theta}) - {h}(\theta) \right)^{<\alpha - 1>}f(\theta)
 - \left({h}(\theta) \right)^{<\alpha - 1>}g(\theta)\right] d\theta = 0, \,
 k = -1,  -2, \dots.
 \end{equation}

It follows from these equations that the spectral characteristic ${h}(\theta)$ of the estimate is determined by the relation
 \begin{equation}\label{sp-char2}
\left( A(e^{i\theta}) - {h}(\theta) \right)^{<\alpha - 1>}f(\theta)
 - \left({h}(\theta) \right)^{<\alpha - 1>}g(\theta) = \overline{C(e^{i\theta})},
 \end{equation}
 $$ C(e^{i\theta})=\sum_{j = 0}^{\infty} c_j e^{i j \theta},$$
 where $c_j$ are unknown coefficients.
These unknown coefficients $c_j$ are determined from the condition
${h}(\theta)\in L^{\alpha}_{-}(f+g)$
which gives us the system of equations

\begin{equation}\label{speq2}
\int_{-\pi}^{\pi} {e^{-i\theta k}\,\, {h}(\theta)} d\theta = 0,\quad k = 0, 1,\dots.
 \end{equation}

 The variance of the optimal estimate of the functional is calculated by the formula

\begin{equation}\label{var2}
\left\|{A} \xi- \hat{A} \xi \right\|_\alpha^\alpha =
\int_{-\pi}^{\pi}
\left| A(e^{i\theta}) - {h}(\theta) \right|^{\alpha}f(\theta)d \theta
 +
 \int_{-\pi}^{\pi} \left|{h}(\theta) \right|^{\alpha}g(\theta)d \theta.
 \end{equation}

We can conclude that the following theorem holds true.

\begin{thm}\label{thm1}
Let $\{\xi_k,k\in\mathbb Z\}$ and $\{\eta_k,k\in\mathbb Z\}$ be mutually independent harmonizable symmetric $\alpha$-stable random sequences which have
absolutely continuous spectral measures and the spectral densities $f(\theta)>0$ and $g(\theta)>0$ satisfying the minimality condition \eqref{eq:minimality}.
The optimal linear estimate $\hat{A} \xi$ of the functional $A \xi~=~\sum_{j = 0}^{\infty} a_j \xi_j$,
that depends on the unknown values $\xi_j,j=0,1,\dots,$ of the sequence  from observations of the sequence $\{\xi_k+\eta_k,k\in\mathbb Z\}$
at points  $k = -1,  -2, \dots,$ is calculated by formula
\eqref{estim2}. The spectral characteristic ${h}(\theta)$ of the estimate is determined by  equation
\eqref{sp-char2}, where the unknown coefficients $c_j$ are determined from the system of equations
\eqref{speq2}.
The variance of the optimal estimate of the functional is calculated by formula \eqref{var2}.
\end{thm}

\subsection{Extrapolation problem. Observations without noise}

Consider the problem of optimal linear estimation of the functional
$$A \xi~=~\sum_{j = 0}^{\infty} a_j \xi_j = \int_{-\pi}^{\pi} A(e^{i\theta})dZ^{\xi}(\theta),\quad A(e^{i\theta})=\sum_{j = 0}^{\infty}a_j e^{ij\theta},
$$
 that depends on the
unknown values of a harmonizable symmetric $\alpha$-stable random sequence $\{\xi_j,j\in\mathbb Z\},$ from observations of the sequence $\xi_k$
at points $k = -1,  -2, \dots$.
Let $\{\xi_k,k\in\mathbb Z\}$ be a harmonizable symmetric $\alpha$-stable random sequence which has absolutely continuous spectral measure and the spectral density $f(\theta)>0$ satisfying the minimality condition
 \begin{equation}\label{minimality1}
 \int_{-\pi}^{\pi} (f(\theta))^{-1/(\alpha-1)}d\theta<\infty.
 \end{equation}
 The optimal estimate $\hat{A} \xi$ of the functional ${A} \xi$ is of the form
 \begin{equation}\label{estim_f}
 \hat{A} \xi = \int_{-\pi}^{\pi} {h}(\theta) dZ^{\xi}(\theta).
 \end{equation}
 The spectral characteristic ${h}(\theta)$ of the optimal linear estimate $\hat{A} \xi$ of the functional is calculated by the formula
 \begin{equation}\label{spchar_f}
 {h}(\theta) = A(e^{i\theta}) - \left(\overline{C(e^{i\theta})}) \right)^{<\frac{1}{\alpha - 1}>} \left(f(\theta)\right)^{\frac{-1}{\alpha - 1}},
 \end{equation}
 where the unknown coefficients $c_j$, $j = 0, 1,\dots$, are determined from the system of equations
 \begin{equation}\label{eq_sphar_f}
 \int_{-\pi}^{\pi} e^{-i\theta k} \left( \left(\sum_{j = 0}^{\infty}a_j e^{i j \theta}\right) - \left(\sum_{j = 0}^{\infty} \overline{c}_j e^{-i j \theta} \right)^{<\frac{1}{\alpha - 1}>} \left( f(\theta) \right)^{\frac{-1}{\alpha - 1}}\right) d\theta = 0, k = 0, 1,\dots
 \end{equation}
 The variance of the optimal estimate of the functional is calculated by the formula
 \begin{equation}\label{variance_f}
\left\| {A}_N \xi- \hat{A} \xi \right\|_\alpha^\alpha =
\int_{-\pi}^{\pi} \left|\left(\overline{C(e^{i\theta}}) \right)^{<\frac{1}{\alpha - 1}>} (f(\theta))^{\frac{-1}{\alpha - 1}} \right|^{\alpha} f(\theta) d \theta.
\end{equation}

As a corollary from Theorem \ref{thm1} the following statement holds true.

\begin{Corollary}\label{cor1}
Let $\{\xi_k,k\in\mathbb Z\}$ be a harmonizable symmetric $\alpha$-stable random sequence which has absolutely continuous spectral measure and the spectral density $f(\theta)>0$ satisfying the minimality condition
\eqref{minimality1}.
The optimal linear estimate $\hat{A} \xi$ of the functional $A \xi~=~\sum_{j = 0}^{\infty} a_j \xi_j$,
that depends on the unknown values $\xi_j,j=0,1,\dots,$ of the sequence  from observations of the sequence $\{\xi_k,k\in\mathbb Z\}$
at points  $k = -1,  -2, \dots,$ is of the form \eqref{estim_f}.
 The spectral characteristic ${h}(\theta)$ of the optimal linear estimate $\hat{A} \xi$ of the functional is calculated by  formula
 \eqref{spchar_f},
 where the unknown coefficients $c_j$, $j = 0, 1,\dots$, are determined from the system of equations
 \eqref{eq_sphar_f}.
 The variance of the optimal estimate of the functional is calculated by  formula
\eqref{variance_f}.
\end{Corollary}

\subsection{Extrapolation problem. Stationary sequences}

Consider the problem in the particular case where $\alpha=2$.
In this case the harmonizable symmetric $\alpha$-stable random sequences $\{\xi_k,k\in\mathbb Z\}$ and $\{\eta_k,k\in\mathbb Z\}$
are stationary sequences and we have the problem of the optimal estimation of the linear functional $A \xi~=~\sum_{j = 0}^{\infty} a_j \xi_j$
that depends on the unknown values $\xi_j,j=0,1,\dots,$ of a stationary random sequence
from observations of the stationary sequence $\{\xi_k+\eta_k,k\in\mathbb Z\}$ at points $k = -1,  -2, \dots$.

We will suppose that stationary sequences $\{\xi_k,k\in\mathbb Z\}$ and $\{\eta_k,k\in\mathbb Z\}$ have
 spectral densities $f(\theta)>0$ and $g(\theta)>0$ satisfying the minimality condition
\begin{equation}\label{eq:minimality}
 \int_{-\pi}^{\pi} (f(\theta)+g(\theta))^{-1}d\theta<\infty.
\end{equation}

The optimal linear estimate $\hat{A} \xi$ of the functional $A\xi$ is of the form
\eqref{estim2}, where
the spectral characteristic ${h}(\theta)$ of the optimal estimate and the variance of the optimal estimate, determined by equations
\eqref{sp-char2}, \eqref{var2}, are of the form
\begin{equation} \label{sphar_a2} \begin{split}
&h(\theta)=\frac{A(e^{i\theta})f(\theta)-C(e^{i\theta})}{f(\theta)+g(\theta)}=\\
&=A(e^{i\theta})-\frac{A(e^{i\theta})g(\theta)+C(e^{i\theta})}{f(\theta)+g(\theta)},
\end{split} \end{equation}
\begin{equation} \label{var_a2} \begin{split}
\Delta(h;f,g)=\left\|{A}_N \xi- \hat{A} \xi \right\|_2^2&=\frac{1}{2\pi} \int\limits_{-\pi}^{\pi}
\frac{\left|A(e^{i\theta})g(\theta)+C(e^{i\theta})
\right|^2}{(f(\theta)+g(\theta))^2}f(\theta)d\theta\\
&+\frac{1}{2\pi} \int\limits_{-\pi}^{\pi}\frac{\left|
A(e^{i\theta})f(\theta)-C(e^{i\theta})
\right|^2}{(f(\theta)+g(\theta))^2}g(\theta)d\theta.\\
\end{split} \end{equation}

The unknown coefficients $c_j$, $j=0,1,\dots,$ are determined from the system of equations
\eqref{speq2} which is of the form in this case
\begin{equation*} \begin{split}
\int\limits_{-\pi}^{\pi}\left(A(e^{i\theta})\frac{f(\theta)}{f(\theta)+g(\theta)}-\frac{C(e^{i\theta})}{f(\theta)+g(\theta)}\right)e^{-ik\theta}d\theta=0, \quad k = 0, 1,\dots
\end{split} \end{equation*}
 From this system of equations we get the following equations
 \begin{equation}\label{3} \begin{split}
\sum\limits_{j=0}^{\infty}a_j\int\limits_{-\pi}^{\pi}\frac{e^{i(j-k)\theta}f(\theta)}{f(\theta)+g(\theta)}d\theta
-\sum\limits_{j=0}^{\infty}c_j\int\limits_{-\pi}^{\pi}\frac{e^{i(j-k)\theta}}{f(\theta)+g(\theta)}d\theta=0,\quad k = 0, 1,\dots
\end{split} \end{equation}

Denote by $\bold{B}$, $\bold{R}$, $\bold{Q}$ operators in the space
$\ell_{2}$, which are determined by matrices with elements
$$B_{k,j}=\frac{1}{2\pi} \int\limits_{-\pi}^{\pi}e^{i(j-k)\theta}\frac{1}{f(\theta)+g(\theta)}d\theta;$$
$$R_{k,j}=\frac{1}{2\pi} \int\limits_{-\pi}^{\pi}e^{i(j-k)\theta}\frac{f(\theta)}{f(\theta)+g(\theta)}d\theta;$$
$$Q_{k,j}=\frac{1}{2\pi} \int\limits_{-\pi}^{\pi}e^{i(j-k)\theta}\frac{f(\theta)g(\theta)}{f(\theta)+g(\theta)}d\theta;$$
$$k,j = 0, 1,2,\dots
$$

With the help of the introduced notations we can write equations (\ref{3}) in the form
\begin{equation*}\begin{split}
\sum\limits_{j=0}^{\infty}R_{k,j}a_j=
\sum\limits_{j=0}^{\infty}B_{k,j}c_j,\quad k = 0, 1,2,\dots.
\end{split} \end{equation*}

These equations can be represented in the matrix-vector form
\begin{equation*}\begin{split}
\bold{R}\bold{a}=\bold{B} \bold{c},
\end{split} \end{equation*}
where $\bold{a}=(a_0,a_1,\dots)$, $\bold{c}=(c_0,c_1,\dots)$. The unknown coefficients $c_j,j=0,1,\dots$ form a solution to this equation and can be represented in the form
$$c_j=(\bold{B}^{-1}\bold{R}\bold{a})_j,$$
where $(\bold{B}^{-1}\bold{R}\bold{a})_j$ is the $j$-th element of the vector $\bold{B}^{-1}\bold{R}\bold{a}$.

Finally, the spectral characteristic and the variance of the optimal estimate  are determined by the formulas (for more details see the books by Moklyachuk (2008), Moklyachuk and Masyutka (2012), Moklyachuk and Golichenko  (2016))

 \begin{equation} \label{sphar_a12} \begin{split}
h(\theta)&=\frac{
A(e^{i\theta})f(\theta)-\sum\limits_{j=0}^{\infty}(\bold{B}^{-1}\bold{R}\bold{a})_je^{ij\theta}
}
{f(\theta)+g(\theta)}=\\
&=A(e^{i\theta})-\frac{
A(e^{i\theta})g(\theta)+\sum\limits_{j=0}^{\infty}(\bold{B}^{-1}\bold{R}\bold{a})_je^{ij\theta}
}
{f(\theta)+g(\theta)},
\end{split} \end{equation}

\begin{equation} \label{var_a12} \begin{split}
\Delta(h;f,g)=\left\|{A} \xi- \hat{A} \xi \right\|_2^2
&=\frac{1}{2\pi} \int\limits_{-\pi}^{\pi}\frac{\left|
A(e^{i\theta})g(\theta)+\sum\limits_{j=0}^{\infty}(\bold{B}^{-1}\bold{R}\bold{a})_je^{ij\theta}
\right|^2}{(f(\theta)+g(\theta))^2}f(\theta)d\theta\\
&+\frac{1}{2\pi} \int\limits_{-\pi}^{\pi}\frac{\left|
A(e^{i\theta})f(\theta)-\sum\limits_{j=0}^{\infty}(\bold{B}^{-1}\bold{R}\bold{a})_je^{ij\theta}
\right|^2}{(f(\theta)+g(\theta))^2}g(\theta)d\theta\\
&=\langle\bold{R}\bold{a},\bold{B}^{-1}\bold{R}\bold{a}\rangle+\langle\bold{Q}\bold{a},\bold{a}\rangle,
\end{split} \end{equation}

So, the following theorem holds true.

\begin{thm}\label{thm2}
Let $\{\xi_k,k\in\mathbb Z\}$ and $\{\eta_k,k\in\mathbb Z\}$ be mutually independent stationary random sequences which have
absolutely continuous spectral measures and the spectral densities $f(\theta)>0$ and $g(\theta)>0$ satisfying the minimality condition \eqref{eq:minimality} with $\alpha=2$.
The optimal linear estimate $\hat{A} \xi$ of the functional $A \xi~=~\sum_{j = 0}^{\infty} a_j \xi_j$,
that depends on the unknown values $\xi_j,j=0,1,\dots,$ of the sequence  from observations of the sequence $\{\xi_k+\eta_k,k\in\mathbb Z\}$
at points  $k = -1,  -2, \dots,$ is calculated by the formula
\eqref{estim2}. The spectral characteristic ${h}(\theta)$ of the estimate is calculated by formula
\eqref{sphar_a12}.
The variance of the optimal estimate of the functional is calculated by formula \eqref{var_a12}.
\end{thm}

\subsection{Extrapolation problem. Stationary sequences. Observations without noise}

Consider the problem of optimal linear estimation of the functional
$$A \xi~=~\sum_{j = 0}^{\infty} a_j \xi_j = \int_{-\pi}^{\pi} A(e^{i\theta})dZ^{\xi}(\theta),\quad A(e^{i\theta})=\sum_{j = 0}^{\infty}a_j e^{ij\theta},
$$
that depends on the
unknown values $\xi_j,j=0,1,\dots,$  of a stationary stochastic sequence  from observations of the sequence $\xi_k$
 at points $k = -1,  -2, \dots$.
 Suppose that the stationary stochastic sequence $\{\xi_k,k\in\mathbb Z\}$ has
 the spectral density $f(\theta)>0$ satisfying the minimality condition
\eqref{minimality1} with $\alpha=2$.
 The optimal linear estimate $\hat{A} \xi$ of the functional $A\xi$
 is of the form \eqref{estim_f}, where the spectral characteristic ${h}(\theta)$ of the optimal linear estimate $\hat{A} \xi$ of the functional $A\xi$ is calculated by the formula
 \begin{equation}\label{sphar_a22}
 {h}(\theta) = A(e^{i\theta}) - \left(\sum\limits_{j=0}^{\infty}(\bold{B}^{-1}\bold{a})_je^{ij\theta} \right) \left(f(\theta)\right)^{-1}.
 \end{equation}
 The variance of the optimal estimate of the functional is calculated by the formula
 \begin{equation}\label{var_a22}
\left\| {A} \xi- \hat{A} \xi \right\|_2^2 =
\frac{1}{2\pi} \int_{-\pi}^{\pi} \left|\sum\limits_{j=0}^{\infty}(\bold{B}^{-1}\bold{a})_je^{ij\theta} \right|^2 f^{-1}(\theta)d\theta=<\bold{B}^{-1}\vec{\bold{a}}, \vec{\bold{a}}>,
\end{equation}
where
$\bold{B}$ is operator in the space $\ell_2$, determined by the matrix with elements
\begin {equation*}
 B(k,j)=\,\frac{1}{2\pi}\int\limits_{-\pi}^{\pi}f^{-1}(\theta)e^{i(j-k)\theta}d\theta, \quad k,j=0,1,\dots.
\end{equation*}

The spectral density $f(\theta)>0$ of the stationary sequence $\{\xi_k,k\in\mathbb Z\}$ satisfies the minimality condition \eqref{minimality1} with $\alpha=2$.
For this reason the function $f^{-1}(\theta)$ admits the factorization
\begin{equation}\label{factor}
\frac{1}{f(\theta)}=\sum\limits_{p=-\infty}^{\infty}b_{p}e^{ip\theta}=\left|\sum\limits_{j=0}^{\infty}\psi_{j}e^{-ij\theta}\right|^2
=\left|\sum\limits_{j=0}^{\infty}\varphi_{j}e^{-ij\theta}\right|^{-2}.
\end{equation}

Denote by $\bold{\Psi}$ and $ \bold{\Phi}$ linear operators in the space $\ell_2$ which are determined by matrices with elements $\bold{\Psi}_{i,j}=\psi_{i-j}$, $\bold{\Phi}_{i,j}=\varphi_{i-j}$, for $0\leq j\leq i$, while $\bold{\Psi}_{i,j}=0$, $\bold{\Phi}_{i,j}=0$, for  $0\leq i< j$.
It can be shown that  $\bold{\Psi}\bold{\Phi}=\bold{\Phi}\bold{\Psi}=I$. The operator  $\bold{B}$ can be represented in the form
$\bold{B}=\bold{\Psi}^{'}\overline{\bold{\Psi}}.$ The operator  $\bold{B}^{-1}$ can be represented in the form
$\bold{B}^{-1}=\overline{\bold{\Phi}}\bold{\Phi}^{'}.$

As a corollary we can represent formula \eqref{var_a22} in the form
\begin{equation}\label{var_a222}
\left\| {A} \xi- \hat{A} \xi \right\|_2^2
=\langle\bold{B}^{-1}\vec{\bold{a}}, \vec{\bold{a}}\rangle
=\langle\overline{\bold{\Phi}}\bold{\Phi}^{'}\vec{\bold{a}}, \vec{\bold{a}}\rangle
=\langle\bold{\Phi}^{'}\vec{\bold{a}}, \bold{\Phi}^{'}\vec{\bold{a}}\rangle
=\langle\bold{A}\vec{\varphi}, \bold{A}\vec{\varphi}\rangle
=\|\bold{A}\vec{\varphi}\|^2,
\end{equation}
where the linear operator $\bold{A}$ in the space $\ell_2$ is determined by matrix with elements
$\bold{A}_{i,j}=a_{i-j}$, $i,j=0,1,\dots$, and the vector $\vec{\varphi}=({\varphi}_0,{\varphi}_1, {\varphi}_2,\dots)$ are determined by elements ${\varphi}_j,j=0,1,\dots$ of the factorization
\eqref{factor}.

So, the following theorem holds true.

\begin{thm}\label{thm3}
Let $\{\xi_k,k\in\mathbb Z\}$ be a stationary stochastic sequence which has absolutely continuous spectral measure and the spectral density $f(\theta)>0$ satisfying the minimality condition \eqref{minimality1} with $\alpha=2$.
The optimal linear estimate $\hat{A} \xi$ of the functional $A \xi~=~\sum_{j = 0}^{\infty} a_j \xi_j$,
that depends on the unknown values $\xi_j,j=0,1,\dots,$ of the sequence $\{\xi_k,k\in\mathbb Z\}$ from observations of the sequence $\{\xi_k,k\in\mathbb Z\}$
at points  $k = -1,  -2, \dots$ is of the form \eqref{estim_f}, where the spectral characteristic ${h}(\theta)$ of the optimal linear estimate $\hat{A} \xi$ of the functional is calculated by the formula \eqref{sphar_a22}.
The variance of the optimal estimate of the functional can be calculated by  formula \eqref{var_a22} as well as by formula \eqref{var_a222}.
\end{thm}


\begin{exm} \label{example1}

Consider the problem of optimal linear estimation of the functional
$$A \xi~=~\sum_{j = 0}^{\infty} a_j \xi_j $$
that depends on the
unknown values $\xi_j,j=0,1,\dots,$  of a stationary stochastic sequence  from observations of the sequence $\xi_k$
 at points $k = -1,  -2, \dots$.
 Suppose that the stationary stochastic sequence $\{\xi_k,k\in\mathbb Z\}$ has
 the spectral density
$$f(\theta) = |1 - \alpha e^{-i\theta}|^{-2}.$$
The function $f^{-1}(\theta)= |1 - \alpha e^{-i\theta}|^{2}$ admits the factorization
$$f^{-1}(\theta) = b_{-1} e^{-i\theta} + b_{0}+b_{1} e^{i\theta}
=\left|\sum\limits_{j=0}^{\infty}\psi_{j}e^{-ij\theta}\right|^2
=\left|\sum\limits_{j=0}^{\infty}\varphi_{j}e^{-ij\theta}\right|^{-2},$$
where $b_{0} = 1+|\alpha|^2$, $b_{-1} = -\alpha$, $b_{1} = -\bar{\alpha}$, $b_p=0, |p|>1$ are Fourier coefficients of the function $f^{-1}(\theta)$;
$\psi_{0} = 1, \psi_{1} = -\alpha,  \psi_j=0, j>1$; $\varphi_{j} = \alpha ^{j},  j\geq 0$;
$b_p=\sum\limits_{k=0}^{\infty}\psi_{k}\overline{\psi}_{k+p},$ $p\geq0$, and $b_{-p}=\overline{b_p}$, $p\geq0$.

The optimal linear estimate $\hat{A} \xi$ of the functional $A\xi$
 is of the form \eqref{estim_f}, where the spectral characteristic ${h}(\theta)$ of the optimal linear estimate $\hat{A} \xi$ of the functional $A\xi$ is calculated by the formula
 \begin{equation}\label{sphar_a2222}
 {h}(\theta) = \left(\sum_{j = 0}^{\infty}a_j e^{ij\theta}\right) - \left(\sum\limits_{j=0}^{\infty}(\bold{B}^{-1}\bold{a})_je^{ij\theta} \right)
\left(b_{-1} e^{-i\theta} + b_0+b_1 e^{i\theta}\right),
 \end{equation}

Making use of the relation
${B}^{-1}=\overline{{\Phi}}{\Phi}^{'}$
we find that the matrix $(B)^{-1}$ is of the form
\begin{equation*}
B^{-1}= \left(\begin{array}{ccccc}
  1 & \alpha & \alpha ^2 & \alpha ^3 &  \ldots \\
  \overline{\alpha} & \overline{\alpha}\alpha+1 & \overline{\alpha}\alpha^{2}+\alpha & \overline{\alpha}\alpha^{3}+\alpha^{2} &  \ldots \\
 \overline{\alpha^2} & \overline{\alpha^2}\alpha+\overline{\alpha} & \overline{\alpha^2}\alpha^2+\overline{\alpha}\alpha+1 & \overline{\alpha^2}\alpha^3+\overline{\alpha}\alpha^2+\alpha &  \ldots \\
 \overline{\alpha^3} & \overline{\alpha^3}\alpha+\overline{\alpha^2} & \overline{\alpha^3}\alpha^2+\overline{\alpha^2}\alpha+\overline{\alpha} & \overline{\alpha^3}\alpha^3+\overline{\alpha^2}\alpha^2++\overline{\alpha}\alpha +1&  \ldots \\
 \ldots
\end{array}\right).
\end{equation*}
Consider the problem under the condition $a_j=0,j\geq3$. In this case the coefficients
$c(j)=\left(\bold{B}^{-1}\vec{\bold{a}}\right)_j, j =0, 1, 2,\ldots$ are as follows
$$
c_0=a_0+a_1\alpha+a_2\alpha^2,
$$
$$
c_1=a_0\overline{\alpha}+a_1(\overline{\alpha}\alpha+1)+a_2(\overline{\alpha}\alpha^{2}+\alpha),
$$
$$
c_2=a_0\overline{\alpha^2}+a_1(\overline{\alpha^2}\alpha+\overline{\alpha})+a_2(\overline{\alpha^2}\alpha^2+\overline{\alpha}\alpha+1),
$$
$$
c_j=a_0\overline{\alpha^j}+a_1(\overline{\alpha^j}\alpha+\overline{\alpha^{j-1}})+a_2(\overline{\alpha^j}\alpha^2+\overline{\alpha^{j-1}}\alpha+\overline{\alpha^{j-2}}),j\geq3.
$$
The spectral characteristic ${h}_2(\theta)$ of the optimal linear estimate $\hat{A}_2 \xi$ of the functional $A_2\xi=a_0+a_1\xi_1+a_2\xi_2$ is calculated by the formula
$$
 {h}_2(\theta) = \left(a_0+a_1 e^{i\theta} + a_2 e^{i2\theta} \right) - \left(\sum\limits_{j=0}^{\infty}c_je^{ij\theta} \right)
\left(b_{-1} e^{-i\theta} + b_0+b_1 e^{i\theta}\right)=
$$
$$
=-c_0b_{-1}e^{-i\theta}=\left(a_0\alpha +a_1\alpha^2+a_2\alpha^3  \right)e^{-i\theta}.
$$

The error of the estimate is calculated by the formula
$$
\left\| {A}_2 \xi- \hat{A}_2 \xi \right\|_2^2=\langle\bold{A}\vec{\varphi}, \bold{A}\vec{\varphi}\rangle
=\langle\bold{B}^{-1}\vec{\bold{a}}, \vec{\bold{a}}\rangle=
$$
$$
=
a_0^2+a_0a_1(\alpha+\overline{\alpha})+a_1^2(1+\alpha^2)+
a_0a_2(\alpha^2+\overline{\alpha^2})+
a_1a_2(\alpha+\overline{\alpha})(1+\alpha^2)+
a_2^2(1+\alpha^2+\alpha^4).
$$

\end{exm}

 \section{Extrapolation problem. Minimax approach}

 The value of the error
 $$\Delta\left(h(f,g);f,g\right):= \left\| \hat{A} \xi- A \xi \right\|_\alpha^\alpha $$
 and
the spectral characteristic $h(f,g):={h}(\theta)$ of the optimal estimate $\hat{A}\xi$ of the functional $A\xi$
can be calculated by the proposed formulas
 only in the case where we know the spectral densities $f(\theta)$ and $g(\theta)$ of the harmonizable symmetric $\alpha$-stable stochastic sequences $\{\xi_k,k\in\mathbb Z\},$  and $\{\eta_k,k\in\mathbb Z\}$. However, usually we do not have exact values of the spectral densities of stochastic sequences, while  we often know a set $D=D_f\times D_g$ of admissible spectral densities. In this case we can apply the minimax-robust method of estimation to the extrapolation problem. This method let us find an estimate that minimizes
the maximum of the errors for all spectral densities from the
given set $D=D_f\times D_g$ of admissible spectral
densities simultaneously (see books by Moklyachuk (2008), Moklyachuk and Masyutka (2012), Moklyachuk and Golichenko (2016) for more details).

 \begin{definition}
 For a given class of spectral densities $D=D_f\times D_g$ the spectral densities $f_0(\theta)\in D_f$, $g_0(\theta)\in D_g$ are called the least favorable in $D=D_f\times D_g$ for the optimal linear estimation of the functional $A \xi$ if the following relation holds true
 $$\Delta\left(f_0,g_0\right)=\Delta\left(h\left(f_0,g_0\right);f_0,g_0\right)=\max\limits_{(f,g)\in D_f\times D_g}\Delta\left(h\left(f,g\right);f,g\right).$$
 \end{definition}

 \begin{definition}
 For a given class of spectral densities $D=D_f\times D_g$ the spectral characteristic $h^0=h(f_0,g_0)$ of the optimal estimate $\hat{A}\xi$ of the functional $A\xi$ is called
 minimax (robust)
 for the optimal linear estimation of the functional $A \xi$ if the following relations hold true
$$h^0\in H_D= \bigcap\limits_{(f,g)\in D_f\times D_g} L^{\alpha}(f+g),$$
$$\min\limits_{h\in H_D}\max\limits_{(f,g)\in D}\Delta\left(h;f,g\right)=\max\limits_{(f,g)\in D}\Delta\left(h^0;f,g\right).$$
 \end{definition}

 The least favorable spectral densities $f_0(\theta)$, $g_0(\theta)$ and the minimax spectral characteristic $h^0=h(f_0,g_0)$ form a saddle point of the function $\Delta \left(h;f,g\right)$ on the set $H_D \times D$. The saddle point inequalities
$$\Delta\left(h;f_0,g_0\right)\geq\Delta\left(h^0;f_0,g_0\right)\geq \Delta\left(h^0;f,g\right) $$ $$ \forall h \in H_D, \forall f \in D_f, \forall g \in D_g$$
holds true if $h^0=h(f_0,g_0)$ and $h(f_0,g_0)\in H_D,$ where $(f_0,g_0)$ is a solution to the constrained optimization problem
\begin{equation} \label{condextr}
\max\limits_{(f,g)\in D_f\times D_g}\Delta\left(h(f_0,g_0);f,g\right)=\Delta\left(h(f_0,g_0);f_0,g_0\right),
\end{equation}
\begin{equation}\label{delta41}\begin{split}
\Delta\left(h(f_0,g_0);f,g\right)&= \left\|{A} \xi- \hat{A} \xi \right\|_\alpha^\alpha\\ & =
\int_{-\pi}^{\pi}
\left| A(e^{i\theta}) - {h}^0(\theta) \right|^{\alpha}f(\theta)d \theta
 +
 \int_{-\pi}^{\pi} \left|{h}^0(\theta) \right|^{\alpha}g(\theta)d \theta.
\end{split}\end{equation}

The constrained optimization problem \eqref{condextr} is equivalent to the unconstrained optimization problem
\begin{equation} \label{8}
\Delta_D(f,g)=-\Delta(h(f_0,g_0);f,g)+\delta(f,g\left|D_f\times D_g\right.)\rightarrow \inf,
\end{equation}
where $\delta(f,g\left|D_f\times D_g\right.)$ is the indicator function of the set $D=D_f\times D_g$.
 Solution $(f_0,g_0)$ to the problem (\ref{8}) is characterized by the condition $0 \in \partial\Delta_D(f_0,g_0),$ where $\partial\Delta_D(f_0,g_0)$ is the subdifferential of the convex functional $\Delta_D(f,g)$ at point $(f_0,g_0)$.
This condition makes it possible to find the least favorable spectral densities in some special classes of spectral densities $D=D_f\times D_g$ (Ioffe and Tikhomirov, 1979; Pshenychnyj, 1971; Rockafellar, 1997; Moklyachuk, 2008b).

Note, that the form of the functional $\Delta(h(f_0,g_0);f,g)$ is convenient for application the method of Lagrange multipliers for
finding solution to the problem \eqref{condextr}.
Making use of the method of Lagrange multipliers and the form of
subdifferentials of the indicator functions
we describe relations that determine least favourable spectral densities in some special classes
of spectral densities

Summing up the derived formulas and the introduced definitions we come to conclusion that the following lemmas hold true

\begin{lemma} \label{lem41}
Let $\{\xi_k,k\in\mathbb Z\}$ and $\{\eta_k,k\in\mathbb Z\}$ be mutually independent harmonizable symmetric $\alpha$-stable random sequences which have
absolutely continuous spectral measures and the spectral densities $f_0(\theta)>0$ and $g_0(\theta)>0$ satisfying the minimality condition \eqref{eq:minimality}.
Let the spectral densities $(f_0,g_0)\in D_f\times D_g$ gives a solution to the constrained optimization problem  \eqref{condextr}.
The spectral densities $(f_0,g_0)$ are the least favorable spectral densities in $D_f\times D_g$ and
$h^0=h(f_0,g_0)$ is the minimax spectral characteristic of the optimal linear estimation
 $\hat{A} \xi$ of the functional $A \xi$,
that depends on the unknown values $\xi_j,j=0,1,\dots,$ of the sequence  from observations of the sequence $\{\xi_k+\eta_k,k\in\mathbb Z\}$
at points  $k = -1,  -2, \dots,$
if $h^0=h(f_0,g_0)\in H_D$.
\end{lemma}

\begin{lemma} \label{lem42}
Let $\{\xi_k,k\in\mathbb Z\}$ be a harmonizable symmetric $\alpha$-stable random sequence which has absolutely continuous spectral measure and the spectral density $f_0(\theta)>0$ satisfying the minimality condition
\eqref{minimality1}.
Let the spectral density $f_0\in D_f$ gives a solution to the constrained optimization problem
\begin{equation} \label{condextr42}
\max\limits_{f\in D_f}\Delta\left(h(f_0);f\right)=\Delta\left(h(f_0);f_0\right),
\end{equation}
\begin{equation} \label{delta42}
\Delta\left(h(f_0);f\right)=\left\| {A} \xi- \hat{A} \xi \right\|_\alpha^\alpha =
\int_{-\pi}^{\pi} \left|\left(C^0(e^{i\theta}) \right)^{<\frac{1}{\alpha - 1}>} (f_0(\theta))^{\frac{-1}{\alpha - 1}} \right|^{\alpha} f(\theta) d \theta.
\end{equation}
The spectral density $f_0$ is the least favorable spectral density in $D_f$ and
$h^0=h(f_0)$ is the minimax spectral characteristic of the optimal linear estimate $\hat{A} \xi$ of the functional $A \xi$,
that depends on the unknown values $\xi_j,j=0,1,\dots,$ of the sequence  from observations of the sequence $\{\xi_k,k\in\mathbb Z\}$
at points $k = -1,  -2, \dots,$
if $h^0=h(f_0)\in H_D$.
\end{lemma}

\begin{lemma} \label{lem43}
Let $\{\xi_k,k\in\mathbb Z\}$ and $\{\eta_k,k\in\mathbb Z\}$ be mutually independent stationary random sequences which have
absolutely continuous spectral measures and the spectral densities $f_0(\theta)>0$ and $g_0(\theta)>0$ satisfying the minimality condition \eqref{eq:minimality} with $\alpha=2$.
Let spectral densities $(f_0,g_0)\in D_f\times D_g$ gives a solution to the constrained optimization problem
\begin{equation} \label{condextr43}
\max\limits_{(f,g)\in D_f\times D_g}\Delta\left(h(f_0,g_0);f,g\right)=\Delta\left(h(f_0,g_0);f_0,g_0\right),
\end{equation}
\begin{equation} \label{delta43} \begin{split}
\Delta\left(h(f_0,g_0);f,g\right)
&=\frac{1}{2\pi}\int\limits_{-\pi}^{\pi}\frac{\left|
A(e^{i\theta})g_0(\theta)+\sum\limits_{j=0}^{\infty}((\bold{B}^0)^{-1}\bold{R}^0\bold{a})_je^{ij\theta}
\right|^2}{(f_0(\theta)+g_0(\theta))^2}f(\theta)d\theta\\
&+\frac{1}{2\pi}\int\limits_{-\pi}^{\pi}\frac{\left|
A(e^{i\theta})f_0(\theta)-\sum\limits_{j=0}^{\infty}((\bold{B}^0)^{-1}\bold{R}^0\bold{a})_je^{ij\theta}
\right|^2}{(f_0(\theta)+g_0(\theta))^2}g(\theta)d\theta.\\
\end{split} \end{equation}
The spectral densities $(f_0,g_0)$ are the least favorable spectral densities in $D_f\times D_g$ and
$h^0=h(f_0,g_0)$ is the minimax spectral characteristic of the optimal linear estimate
 $\hat{A} \xi$ of the functional $A \xi$,
that depends on the unknown values $\xi_j,j=0,1,\dots,$ of the sequence  from observations of the sequence $\{\xi_k+\eta_k,k\in\mathbb Z\}$
at points  $k = -1,  -2, \dots,$
if $h^0=h(f_0,g_0)\in H_D$.
\end{lemma}

\begin{lemma} \label{lem44}
Let $\{\xi_k,k\in\mathbb Z\}$ be a stationary random sequence which has
absolutely continuous spectral measure and the spectral density $f_0(\theta)>0$ satisfying the minimality condition \eqref{minimality1} with $\alpha=2$.
Let the spectral density $f_0\in D_f$ gives a solution to the constrained optimization problem
\begin{equation} \label{condextr44}
\max\limits_{f\in D_f}\Delta\left(h(f_0);f\right)=\Delta\left(h(f_0);f_0\right),
\end{equation}
 \begin{equation}\label{delta44}
\Delta\left(h(f_0);f\right) =
\frac{1}{2\pi}\int_{-\pi}^{\pi} \left|\sum\limits_{j=0}^{\infty}((\bold{B}^0)^{-1}\bold{a})_je^{ij\theta} \right|^2 f_0^{-2}(\theta)f(\theta)d\theta.
\end{equation}
The spectral density $f_0$ is the least favorable spectral density in $D_f$ and
$h^0=h(f_0)$ is the minimax spectral characteristic of the optimal linear estimate $\hat{A} \xi$ of the functional $A \xi$,
that depends on the unknown values $\xi_j,j=0,1,\dots,$ of the sequence  from observations of the sequence $\{\xi_k,k\in\mathbb Z\}$
at points  $k = -1,  -2, \dots,$
if $h^0=h(f_0)\in H_D$.
\end{lemma}

 \section{Least favorable spectral densities in the class $D_f^{\beta}\times D_g^{\varepsilon}$}

 Consider the problem of optimal estimation of the linear functional
 $A \xi$ that depends on the unknown values $\xi_j,j=0,1,\dots,$ of a random sequence  from observations of the sequence $\{\xi_k+\eta_k,k\in\mathbb Z\}$
 at points  $k = -1,  -2, \dots$, where $\{\xi_k,k\in\mathbb Z\}$ and $\{\eta_k,k\in\mathbb Z\}$
are mutually independent harmonizable symmetric $\alpha$-stable random sequences which have
 spectral densities $f(\theta)>0$ and $g(\theta)>0$ satisfying the minimality condition \eqref{eq:minimality} from
 the class of admissible spectral densities $D=D_f^{\beta}\times D_g^{v,u}$, where
$$
D_f^{\beta} = \left\{f(\theta)\left|\int\limits_{-\pi}^{\pi} (f(\theta))^{\beta}d\theta= P_1\right.\right\},$$
$$D_g^{\varepsilon} = \left\{g(\theta)\left| g(\theta)=(1-{\varepsilon}) g_1(\theta)+ {\varepsilon}w(\theta),\int\limits_{-\pi}^{\pi} g(\theta)d\theta= P_2\right.\right\}. $$
Assume that spectral densities $f_0\in D_f^{\beta}$,
$g_0\in D_g^{\varepsilon}$ and the functions $h_f(f_0,g_0)$, $h_g(f_0,g_0)$, determined by the equations
\begin{equation} \label{hf51}
\qquad\qquad h_f(f_0,g_0)=\left| A(e^{i\theta}) - {h}^0(\theta) \right|^{\alpha},
\end{equation}
\begin{equation} \label{hg51}
h_g(f_0,g_0)=\left|{h}^0(\theta) \right|^{\alpha},
\end{equation}
 \begin{equation}\label{sp51}
\left( A(e^{i\theta}) - {h}^0(\theta) \right)^{<\alpha - 1>}f_0(\theta)
 - \left({h}^0(\theta) \right)^{<\alpha - 1>}g_0(\theta) = C^0(e^{i\theta}),
 \end{equation}
\begin{equation}\label{speq51}
\int_{-\pi}^{\pi} {e^{-i\theta k}\,\, {h}^0(\theta)} d\theta = 0,\quad k = 0, 1,\dots
 \end{equation}
are bounded.
Under these conditions the functional
\begin{equation}\begin{split}
\Delta\left(h(f_0,g_0);f,g\right)&= \left\|{A} \xi- \hat{A} \xi \right\|_\alpha^\alpha\\ & =
\int_{-\pi}^{\pi}
\left| A(e^{i\theta}) - {h}^0(\theta) \right|^{\alpha}f(\theta)d \theta
 +
 \int_{-\pi}^{\pi} \left|{h}^0(\theta) \right|^{\alpha}g(\theta)d \theta.
\end{split}\end{equation}
is linear and continuous in the $L_1\times L_1$ space
and we can apply the Lagrange multipliers method to derive that
the least favorable
densities $f_0\in D_f^{\beta}$, $g_0\in D_g^{\varepsilon}$ satisfy the equations
\begin{equation} \label{eqf51}
\left| A(e^{i\theta}) - {h}^0(\theta) \right|^{\alpha}= \gamma_1 \left(f_0(\theta)\right)^{\beta - 1},
\end{equation}
\begin{equation} \label{eqg51}
\left|{h}^0(\theta) \right|^{\alpha}= \left(\varphi_1 (\theta)+\gamma_2 \right),
\end{equation}
where
$\varphi_1 (\theta)\leq 0$ and $\varphi_1(\theta)= 0$ if $g_0(\theta)\geq (1-{\varepsilon}) g_1(\theta)$;
$\gamma_1,$ $\gamma_2$ are the Lagrange multipliers which are determined from the conditions
$$
\int\limits_{-\pi}^{\pi} (f_0(\theta))^{\beta}d\theta= P_1,\quad \int\limits_{-\pi}^{\pi} g_0(\theta)d\theta= P_2.
$$

Thus, the following statement holds true.

\begin{thm}\label{thm51} Let the spectral densities $f_0\in D_f^{\beta}$, $g_0\in D_g^{\varepsilon}$ satisfy the minimality condition \eqref{eq:minimality} and let the functions
$h_f(f_0,g_0)$, $h_g(f_0,g_0)$ determined by formulas \eqref{hf51}, \eqref{hg51}, \eqref{sp51}, \eqref{speq51}
be bounded.
The spectral densities
$f_0(\theta)$ and $g_0(\theta)$ are the least favorable in
the class $D=D_f^{\beta}\times D_g^{\varepsilon}$ for the optimal linear
estimation of the functional $A \xi$ if they satisfy equations
\eqref{eqf51}, \eqref{eqg51}
and determine a solution to the optimization problem \eqref{condextr}.
The minimax-robust spectral characteristic $h(f_0,g_0)$ of the optimal estimate of the
functional $A \xi$ is determined by formulas \eqref{sp51}, \eqref{speq51}.
\end{thm}

 \subsection{Least favorable spectral densities. Observations without noise}

Consider the problem of optimal linear estimation of the functional
 $A \xi= \sum_{j = 0}^{\infty} a_j \xi_j$ that depends on the
unknown values $\xi_j, j = 0, 1, \ldots$, from observations of the sequence $\{\xi_k,k\in\mathbb Z\}$
 at points $k = -1,  -2, \dots$, where $\{\xi_k,k\in\mathbb Z\}$ is a harmonizable symmetric $\alpha$-stable random sequence which have
spectral density $f_0(\theta)>0$ satisfying the minimality condition \eqref{minimality1} from
 the class of admissible spectral densities $D_f^{\beta}$.
Assume that spectral density $f_0\in D_f^{\beta}$ and the function $h_f(f_0)$ determined by the equation
\begin{equation} \label{hf52}
h_f(f_0)=
\left|\left(C^0(e^{i\theta}) \right)^{<\frac{1}{\alpha - 1}>} (f_0(\theta))^{\frac{-1}{\alpha - 1}} \right|^{\alpha}
\end{equation}
is bounded.
Under this condition the functional
\eqref{delta42} is linear and continuous in the $L_1$ space
and we can apply the  method of Lagrange multipliers to
find solution of the
constrained optimization problem
\eqref{condextr42}
and derive that
the least favorable
density $f_0\in D_f^{\beta}$ satisfy the equation
\begin{equation} \label{eqf52}
\left|\left(C^0(e^{i\theta}) \right)^{<\frac{1}{\alpha - 1}>} (f_0(\theta))^{\frac{-1}{\alpha - 1}} \right|^{\alpha}= \gamma_1 \left(f_0(\theta)\right)^{\beta - 1},
\end{equation}
where $\gamma_1$ is the Lagrange multipliers. From this equation we find that the least favorable
density is of the form
\begin{equation}\label{eqf522}
f_0(\theta) = C\left|\sum_{j = 0}^{\infty} c_j e^{-i j \theta} \right|^\frac{\alpha}{\alpha + (\alpha-1)(\beta - 1)}.
\end{equation}
The unknown constants are determined from the optimization  problem \eqref{condextr42} and from the condition
$$
\int\limits_{-\pi}^{\pi} (f_0(\theta))^{\beta}d\theta= P_1.
$$
In the case $\beta=1$ the least favorable
density is of the form
\begin{equation}\label{eqf5222}
f_0(\theta) = C\left|\sum_{j = 0}^{\infty} c_j e^{-i j \theta} \right|.
\end{equation}

The following statement holds true.

\begin{thm}\label{thm52} Let the spectral density $f_0\in D_f^{\beta}$ satisfy the minimality condition \eqref{minimality1} and let the function
$h_f(f_0)$ determined by formula \eqref{hf52}
be bounded.
The spectral density $f_0(\theta)$ is the least favorable in
the class $D_f^{\beta}$ for the optimal linear
estimation of the functional $A \xi$ if it is of the form
\eqref{eqf522}
and determine a solution to the optimization problem \eqref{condextr42}.
The minimax-robust spectral characteristic $h(f_0)$ of the optimal estimate of the
functional $A \xi$ is determined by formula \eqref{spchar_f}.
\end{thm}

\subsection{Least favorable spectral densities. Stationary sequences}

 Consider the problem of the optimal estimation of the linear functional
 $A \xi$ that depends on the unknown values $\xi_j,j=0,1,\dots,$ of a random sequence $\{\xi_k,k\in\mathbb Z\}$  from observations of the sequence $\{\xi_k+\eta_k,k\in\mathbb Z\}$
 at points of time $k = -1,  -2, \dots$, where $\{\xi_k,k\in\mathbb Z\}$ and $\{\eta_k,k\in\mathbb Z\}$
are mutually independent stationary stochastic sequences which have
 spectral densities $f(\theta)>0$ and $g(\theta)>0$ satisfying the minimality condition \eqref{eq:minimality} with $\alpha=2$ from
 the class of admissible spectral densities $D=D_f^{\beta}\times D_g^{\varepsilon}$.

Assume that spectral densities $f_0\in D_f^{\beta}$,
$g_0\in D_g^{\varepsilon}$ and the functions $h_f(f_0,g_0)$, $h_g(f_0,g_0)$, determined by the equations
\begin{equation} \label{hf53}
 h_f(f_0,g_0)=\frac{\left|A(e^{i\theta})g_0(\theta)+\sum\limits_{j=0}^{\infty}((\bold{B}^0)^{-1}\bold{R}^0\bold{a})_je^{ij\theta}\right|^2}{(f_0(\theta)+g_0(\theta))^2},
\end{equation}
\begin{equation} \label{hg53}
h_g(f_0,g_0)=\frac{\left|A(e^{i\theta})f_0(\theta)-\sum\limits_{j=0}^{\infty}((\bold{B}^0)^{-1}\bold{R}^0\bold{a})_je^{ij\theta}
\right|^2}{(f_0(\theta)+g_0(\theta))^2}
\end{equation}
 are bounded.
 Under these conditions the functional
\eqref{delta43} is linear and continuous in the $L_1\times L_1$ space
and we can apply the Lagrange multipliers method to
find solution of the
constrained optimization problem
\eqref{condextr43}
and derive that
the least favorable
densities $f_0\in D_f^{\beta}$, $g_0\in D_g^{\varepsilon}$ satisfy the equations
\begin{equation} \label{eqf53}
{\left|A(e^{i\theta})g_0(\theta)+\sum\limits_{j=0}^{\infty}((\bold{B}^0)^{-1}\bold{R}^0\bold{a})_je^{ij\theta}\right|^2}=\gamma_1
{(f_0(\theta)+g_0(\theta))^2} \left(f_0(\theta)\right)^{\beta - 1},
\end{equation}
\begin{equation} \label{eqg53}
{\left|A(e^{i\theta})f_0(\theta)-\sum\limits_{j=0}^{\infty}((\bold{B}^0)^{-1}\bold{R}^0\bold{a})_je^{ij\theta}
\right|^2}=
{(f_0(\theta)+g_0(\theta))^2} \left(\varphi_1 (\theta)+\gamma_2 \right),
\end{equation}
where
$\varphi_1 (\theta)\leq 0$ and $\varphi_1(\theta)= 0$ if $g_0(\theta)\geq (1-{\varepsilon}) g_1(\theta)$;
$\gamma_1,$ $\gamma_2$ are the Lagrange multipliers which are determined from the conditions
$$
\int\limits_{-\pi}^{\pi} (f_0(\theta))^{\beta}d\theta= P_1,\quad \int\limits_{-\pi}^{\pi} g_0(\theta)d\theta= P_2.
$$

Thus, the following statement holds true.

\begin{thm}\label{thm53}
Let the spectral densities $f_0\in D_f^{\beta}$, $g_0\in D_g^{\varepsilon}$ satisfy the minimality condition \eqref{eq:minimality} with $\alpha=2$ and let the functions
$h_f(f_0,g_0)$, $h_g(f_0,g_0)$ determined by formulas \eqref{hf53}, \eqref{hg53} be bounded.
The spectral densities
$f_0(\theta)$ and $g_0(\theta)$ are the least favorable in
the class $D=D_f^{\beta}\times D_g^{\varepsilon}$ for the optimal linear
estimation of the functional $A \xi$ if they satisfy equations
\eqref{eqf53}, \eqref{eqg53}
and determine a solution to the constrained optimization problem \eqref{condextr43}.
The minimax-robust spectral characteristic $h(f_0,g_0)$ of the optimal estimate of the
functional $A \xi$ is determined by formulas \eqref{sphar_a12}.
\end{thm}

\subsection{Least favorable spectral densities. Stationary sequences. Observations without noise}

 Consider the problem of the optimal linear estimation of the functional
 $A \xi$ that depends on the
unknown values of a random sequence $\{\xi_j,j\in\mathbb Z\}$ from observations of the sequence
 at points $k = -1,  -2, \dots$, where the stationary random sequence $\{\xi_k,k\in\mathbb Z\},$ has
 the spectral density $f(\theta)>0$ satisfying the minimality condition \eqref{minimality1} with $\alpha=2$
 from the class of admissible spectral densities $D_f^{\beta}$.
 Assume that spectral density $f_0\in D_f^{\beta}$ and the function $h_f(f_0)$ determined by the equation
\begin{equation} \label{hf54}
 h_f(f_0)=\left|\sum\limits_{j=0}^{\infty}((\bold{B}^0)^{-1}\bold{a})_je^{ij\theta} \right|^2f_0^{-2}(\theta)
\end{equation}
is bounded.

Under this condition the functional
\eqref{delta44} is linear and continuous in the $L_1$ space
and we can apply the Lagrange multipliers method to
find solution of the
constrained optimization problem
\eqref{condextr44}
and derive that
the least favorable
density $f_0\in D_f^{\beta}$ satisfy the equation
\begin{equation} \label{eqf54}
\left|\sum\limits_{j=0}^{\infty}((\bold{B}^0)^{-1}\bold{a})_je^{ij\theta} \right|^2f_0^{-2}(\theta)= \gamma_1 \left(f_0(\theta)\right)^{\beta - 1},
\end{equation}
where $\gamma_1$ is the Lagrange multiplier which is determined from the conditions
$$
\int\limits_{-\pi}^{\pi} (f_0(\theta))^{\beta}d\theta= P_1.
$$

Thus, the following statement holds true.

\begin{thm}\label{thm54}
Let the spectral density $f_0\in D_f^{\beta}$ satisfy the minimality condition \eqref{minimality1} with $\alpha=2$ and let the function
$h_f(f_0)$ determined by formula \eqref{hf54} be bounded.
The spectral density
$f_0(\theta)$ is the least favorable in
the class $D_f^{\beta}$ for the optimal linear
estimation of the functional $A \xi$ if it satisfies equation
\eqref{eqf54}
and determine a solution to the constrained optimization problem \eqref{condextr44}.
The minimax-robust spectral characteristic $h(f_0)$ of the optimal estimate of the
functional $A \xi$ is determined by formulas \eqref{sphar_a22}.
\end{thm}

In the case of $\beta=1$  the set of admissible spectral densities $D=D_f$ is of the form
$$
D_f = \left\{f(\theta)\left|
\frac{1}{2\pi}
\int\limits_{-\pi}^{\pi} f(\theta)d\theta= P\right.\right\}.
$$
Stationary sequences with the spectral densities from such $D_f$ have finite dispersion
$E|\xi_j|^2= P$ and can be represented as a sum of a regular sequence and a singular sequence.
The least favorable in $D_f$ spectral density is density of the regular sequence since singular sequences have zero value of the mean square error of extrapolation.
Spectral densities from $D_f$ of the regular sequences admit the factorization
\begin{equation}\label{factor2222}
f(\theta) = \left|\sum\limits_{j=0}^{\infty}\varphi_{j}e^{-ij\theta}\right|^{2},\quad
\sum\limits_{j=0}^{\infty}|\varphi_{j}|^{2}= P,
\end{equation}
and we can use the following optimization problem to find the least favourable spectral density in the set $D_f$
\begin{equation}\label{var_a2222}
\|\bold{A}\vec{\varphi}\|^2\to\max,\quad \|\vec{\varphi}\|^2=
\sum\limits_{j=0}^{\infty}|\varphi_{j}|^{2}= P,
\end{equation}
where the linear operator $\bold{A}$ in the space $\ell_2$ is determined by matrix with elements
$\bold{A}_{i,j}=a_{i-j}$, $i,j=0,1,\dots$, and the vector $\vec{\varphi}=({\varphi}_0,{\varphi}_1, {\varphi}_2,\dots)$ are determined by elements ${\varphi}_j,j=0,1,\dots$ of the factorization
\eqref{factor2222}.
Solution to the optimization problem \eqref{var_a2222} gives the eigenvector $\vec{\varphi}^0=({\varphi}_0^0,{\varphi}_1^0, {\varphi}_2^0,\dots)$ which corresponds to greatest eigenvalue
$\nu^0$ of the linear operator $\bold{A}$.

We present this result as a theorem

\begin{thm}\label{thm544}
Spectral density
$f_0(\theta)$ is the least favourable in
the class $D_f$ for the optimal linear
estimation of the functional $A \xi$ if it is of the form
\eqref{factor2222}, where $\vec{\varphi}^0=({\varphi}_0^0,{\varphi}_1^0, {\varphi}_2^0,\dots)$ is the eigenvector of the linear operator $\bold{A}$
which corresponds to greatest eigenvalue $\nu^0$ of the operator.
The optimal minimax linear estimate $ \hat{A} { \xi}$ of the functional $A {\xi}=\sum_{j=0}^{ \infty} {a}_j\xi_j$ is of the form
 \[ \hat{A} { \xi}= \sum_{j=0}^{ \infty} {a}_j \left[ \sum_{u=- \infty}^{-1} \varphi_{j-u}^0 { \varepsilon}_u \right] , \]
\noindent where $ \varepsilon_u$
is a standard stationary sequence with orthogonal
values ("white noise" sequence), the sequence $ \{\varphi_u^0: u= 0,1,\dots \}$ is uniquely determined by
coordinates of the eigenvector of the operator $\bold{A}$ that corresponds
to the greatest eigenvalue $\nu^0$ and condition $E \left \|
{ \xi}_j \right \| ^{2} =P$.
\end{thm}

 \begin{exm} \label{exm2.1}
Consider the problem of optimal linear stimulation of the functional
 \[A_{1} { \xi}= \xi (0)+  2\xi (1) \]
that depends on the unknown values $ { \xi}(0)$, $ { \xi}(1)$ of a stationary sequence $ { \xi}(j)$, that satisfies conditions
 \[ E{ \xi}(j)=0,\quad E\left | { \xi}(j) \right |^{2} \le P, \]
based on observations of the sequence $ { \xi}(j)$ at points $j=-1,-2,\dots$.

\noindent Eigenvalues of the operator $\bold{A}_1$ are equal to $(1\pm \sqrt{17})/2 $.
So the greatest eigenvalue is $ \nu_{1} =(1 + \sqrt{17})/2 $. The eigenvector corresponding to the eigenvalue $ \nu_{1} =(1+ \sqrt{17})/2  $ is of the
form $ \vec\varphi = \left \{ \varphi (0), \varphi(1) \right \}$, where
 \[ \varphi (0)= \sqrt{1/2+ 1/2\sqrt{17}},\quad \varphi (1)= \sqrt{1/2- 1/2\sqrt{17}} . \]

The least favourable in
the class $D_f$ for the optimal linear
estimation of the functional $A_1 \xi$
spectral density
is of the form
$$
f(\theta) = P\left|\varphi_{0}+\varphi_{1}e^{-ij\theta}\right|^{2},\quad
\left|\varphi_{0}\right|^{2}+\left|\varphi_{1}\right|^{2}= 1,
$$
which is the spectral density of the least favourable stationary sequence $ { \xi}(j)$ that is a moving average sequence of the form
 \[ { \xi}(j)= \sqrt{P}\varphi (0) \varepsilon (j)+ \sqrt{P}\varphi (1) \varepsilon (j-1)= \]
 \[= \sqrt{P}\sqrt{1/2+ 1/2\sqrt{17}} \,\,\varepsilon(j) +
 \sqrt{P}\sqrt{1/2- 1/2\sqrt{17}}\,\,\varepsilon(j-1), \]
 $\varepsilon(j)$ is a ''white noise'' sequence.
\noindent The optimal linear minimax estimate $ \hat{A}_{1} { \xi}$ of the functional $A_{1} { \xi}$ is of the form
 \[ \hat{A}_{1} { \xi}=  \sqrt{P}\varphi (1)\,\, \varepsilon (-1)= \sqrt{P}\sqrt{1/2- 1/2\sqrt{17}}\,\,\varepsilon (-1). \]

\noindent The mean-square error of the optimal estimate of the functional
$A_{1} { \xi}$ does not exceed $(9 + \sqrt{17})/2$.
 \end{exm}

 \section{Conclusion}
 We propose methods of solution the optimal linear estimation problem for the linear functional
 $A \xi~=~\sum_{j = 0}^{\infty} a_j \xi_j$ that depends on the
unknown values of a random sequence $\{\xi_j,j\in\mathbb Z\}$ from observations of the sequence
 $\{\xi_k+\eta_k,k\in\mathbb Z\}$ at points
$k = -1,  -2, \dots$, where $\{\xi_k,k\in\mathbb Z\}$ and $\{\eta_k,k\in\mathbb Z\}$ are mutually independent harmonizable symmetric $\alpha$-stable random sequences
which have the spectral densities $f(\theta)>0$ and $g(\theta)>0$ satisfying the minimality condition.
The problem is investigated under the condition of spectral certainty as well as under the condition of spectral uncertainty.
 Formulas for calculation the value of the error and spectral characteristic of the optimal linear
estimate of the functional are derived under the condition of spectral certainty where spectral densities of the sequences are exactly known.
 In the case where spectral densities of
the sequences are not exactly known while a set of admissible spectral densities is available, relations which determine least favorable densities and the minimax-robust spectral characteristics  are found for different classes of spectral densities.

\end{document}